# MATHEMATICS: ABSTRACTION AND REALITY
## A SKETCH TOWARD DEEPER ANALYSIS[1]

Radoslav Dimitric

While relaxing at the bank of a lake and enjoying watching ducks that were happily swimming around, a question came to me as to whether a duck with ten ducklings would notice if one of them was quietly removed from the flock, and if not, whether it would notice if six ducklings were removed, etc... Would this represent only a change in mass for the duck's vision, or would this be a completely new quality namely that of counting? I soon realized that the related experiments with animals I had envisaged have been done before and are still being done. Here is one of them: One pigeon was fed a grain of corn in different time intervals, but every seventh kernel was glued to the plate. After some time, the pigeon refused to take every seventh kernel it received. Jackdaws (Corvus monedula) were able to learn what lids to remove from containers each containing up to two food items, until they obtained a given number of rewards, although the trials had different distributions of food by container (Köhler [8]). One of the interpretations was that this is a specific counting without numbers. Similarly the chimpanzees which were "counting" straws could manage to keep exact track only to the fifth number, but some other experiments have them go up to seven (Hinde & Stevenson-Hinde [6]). Do animals have the ability for abstraction necessary for counting? Number seven seems to be the greatest number the animals can master!

Hundreds of thousands of years must have gone by in a distant past, before a man (a manlike ape) understood that nine stones, nine fish and nine trees represent in some sense one entity, in this case united by number nine. A giant qualitative leap in human thinking happened with the introduction of the first mathematical abstractions, with the creation of the Homo Mathematicus who was gifted with a special apparatus for the abstraction of reality as well as its own thoughts as a part of that reality (we view the psyche here as reality, for there is nothing more real than a hallucination or a, mirage, if the latter is not seen exclusively as an optical phenomenon). This ancient feeling of the mankind for mathematics may be witnessed in the societies of Australian Aborigines, in the existence of such a simple, and at the same time perfect, tool like boomerang. Van der Waerden claims in one of his recent studies [16] that the Pythagoras' theorem was known in the ice ages. Plimpton 322, a Babylonian clay tablet ca 1800 years B.C. (written in cuneiform script and kept at Columbia University Rare Book and Manuscript Library) contains 15 Pythagorean triples.

---

1  This is a translation of a paper delivered in Mathematics History and Philosophy Seminar at the Mathematics Institute, Belgrade, 1983/84. Subsequently a number of monographs appeared with related and more elaborate examination of the issues that we treat here. We mention Mac Lane's monograph [9] and Dieudonné's [2].

Going through different stages of abstraction, in its thousands of years of development mathematics has gone through the problem of notation and recording of the natural numbers (the oldest written mathematical document dating back to ca 1850 before Christ is the Moscow mathematical papyrus kept in the Pushkin State Museum of Fine Arts (a larger papyrus, about 1650 BCE is the "Calculating book of the scribe Ahmes;" it is kept in the British Museum). One of the ways the Mayas used to record the numbers were as follows (with the Hindu-Arabic equivalents we use today) (Stingl, [15]):

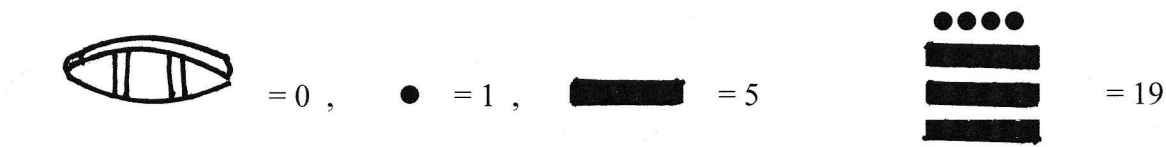

It seems clear that the animals do not have a genetic base that would enable them to mathematicise reality or apply the required abstractions. But why is it that primitive tribes or small children do not possess developed mathematics, although they are biologically equipped with the apparatus able for mathematical thinking equal to the apparatus of the contemporary mathematics expert? The answer to this question is to be found in arguments that follow.

### Culture, or is mathematical reality outside or within us?

In a long history of mathematics there existed supporters of both of the opinions in the question of the title. We quote a distinguished mathematician G. H. Hardy [5]: "I believe that mathematical reality lies outside us, that our function is to discover or observe it, and that the theorems which we prove, and which we describe grandiloquently as our "creations," are simply the notes of our observations." The opposite opinions on scientific theories as purely mental constructions are equally represented, to mention only one representative of this opinion - a mathematical physicist Albert Einstein. Erwin Schrödinger had this to say (from his inaugural lecture, Dec 9, 1922, see [12]) about the role of traditions (read: culture) in forming a scientific theory: "Whence arises the widespread belief that the behavior of molecules is determined by absolute causality, whence the conviction that the contrary is *unthinkable*? Simply through the *custom* inherited through thousands of years of *thinking causally*, which makes the idea of undetermined events, of absolute, primary causality, seem complete nonsense, logical absurdity."

After all, we may ask whether the question of mathematical reality being outside or within us, is well posed at all! The totality of mathematics is certainly a part of the human culture and as such, it is outside of a man, but it is within himself also. All people inherit the culture of their predecessors and every man learns the knowledge from his culture, resulting in the fact that different cultures have different mathematics; for instance, some cultures count in tens, some in twelves (dozens), and some have no names for numbers greater than five. Hottentots had their own mathematics, just as England of Newton's time did. But had Newton had grown up among the Hottentots, he would have had certainly thought like a Hottentot and would have had like mathematical results. We interpret extremely lightly the consequences of our cultural habits, as simple and unique reflections of reality. For example, the Western white men tend to think that the yellow, blue and green pigments are trivial characteristics of the outside world, until we learn that the Creek and Natchez Indians did not distinguish yellow and green, namely they had only one name for both. Also Choctaw, Tunica and Keres among the Pueblo

Indians and a number of other peoples do not make a terminological distinction between green and blue. Do then these colors exist for them?

It is the culture factor that may be used to explain the phenomenon that the same mathematical discovery is arrived at independently, in several different places (several different minds) simultaneously, thus under the influence of very similar cultural climate.

New syntheses in mathematics are formed through the interaction of the ideas in the human nervous system, under the influence of conditions of mathematical culture in a given environment. Culture is the answer to the question why primitive tribes and children do not have developed mathematical minds. Also, no matter how well we teach any single animal to count, its successors will not assimilate counting, because the animals do not have culture either in the form of written or oral traditions.

Influence of any individual brain is a very important factor and cannot be undermined. Just as Homo sapiens is equipped with a superior thinking apparatus in comparison to other animal species, the individuals follow the same pattern - just as some have better sense of smell or hearing, the same can be said about the gray-white matter.

However, several thousands of years before Leibnitz and Newton (in evolutionary terms this is a small period for the brain to undergo changes in its genetics and biology) there existed brains structurally equally valuable as those of Leibnitz' and Newton's, but the lightening appeared in Leibnitz' and Newton's heads, because they were in a cultural environment that enabled this to happen. Had they herded sheep instead of working on science (we are lucky they did not), some other equally exceptional instrument (brain) would be found in which the theory would light up.

We think that the great minds are uniformly distributed under all skies, color of skin, cultural and political traditions (the strength of those minds is most likely distributed according to the Gaussian normal probability law), which means that the problem is not how to make great individuals, but how to make great societies. Great societies by themselves churn out men for great circumstances.

The Germans have arrived, in the past, at the methods that led to deeper trace in the mine of science, giving up sporadic methods of learning. In their technical schools and universities the advance did not have to wait a random genius or a random lucky thought. The whole world admires their learning level in the nineteen century. This discipline of knowledge can be applied outside technology to the "pure" science, and outside science - to the general learning. It constitutes exactly that passage from an amateur to a professional, so necessary in science and especially in mathematics. There have always existed people who devoted their lives to certain areas of thinking. However, the full self-conscious understanding of the power of professionalism in many chapters of knowledge and ways for generating experts was realized for the first time in the nineteen century and, among other countries, mainly in Germany.

Mathematics has had in its history (and will undoubtedly continue to have) several different paths of development. One of the ways of mathematical living is a consistent following of the paths within mathematics theory, without interactions with the "outside" world. We will now elaborate on this very path.

# Mathematics as an abstraction

The science of pure mathematics can claim the right to be the most original creation of human mind, occupying a unique rank among all the creations. Music can claim that right, right after her. In its full depth, mathematics is accessible to the fewest and it is a true art just like music. Both are turned inwards, both draw inspiration from solitude, with a dint of hard work, both follow formal conventions in their development.

The primacy of mathematics comes from the fact that it shows the extremely non-obvious relationships among things. Thus the ideas that are present in the conscious minds of contemporary mathematicians are rather distant from any representations that can be induced by a sensory perception, if that perception is not accompanied in advance by a mathematical knowledge. A great number of questions that keep mathematicians of the present time busy (like the squaring of the circle in old Greece), such as Riemann hypothesis or problems of modern algebra, would likely appear to the ancient mathematicians (who sought simple concrete numbers for their results), as a humorous, somewhat queer games. And for the wide, popular circles this is undoubtedly still the case. There is nothing less popular than modern mathematics with abstraction as its main tool. This also represents something of the symbol of infinite distance. Long time had passed from "innocent games" of George Boole with logical symbols, when his "Laws of Thought" [1] had induced mild teasing of his contemporaries who considered it all to be useless scholar's thinking, to todays insatiable hysteria of the laity and its unlimited trust in computers and everything connected to them. It seems that the magic number of hundred years is almost necessary as a distance from the moment of an abstract treatment of a mathematical problem, to its economic exploitation.

Essential for mathematics is the freeing in it from the special case, so that no mathematical truths are only valid for fish, or for stones or only for people. As long as one is doing mathematics, one is in the area of absolute abstraction. One only claims that if some abstract conditions are satisfied, then some other abstract conditions are satisfied.

The gift of generalization is probably one of the greatest of all gifts. To see particular and seemingly disparate phenomena as a consequence of unifying and clear facts, is something very important and also something to which not, enough attention and credit are paid. The tendency today is to explore very particular notions, but only the unification of those notions shines a new light on the matter, and opens significantly new horizons (Dimitric [3]). Pythagoras, who lived in the sixth century B.C. And whom we do not know much about, was the first to emphasize the importance of generality in reasoning and he sensed the importance of numbers as a tool for building any understanding of the order of nature. The foundation of Pythagorean brotherhood and mystical rumors about their rituals and his influence, offer some proofs that Pythagoras had sensed possible importance of mathematics for the formation of science.

The most impressive fact is that the more mathematics gradually retreats into the upper regions of ever more abstract thought, the more and more, when returning to the earth, it becomes important for the analysis of concrete facts. The civilization unable to overcome its current abstractions, is doomed to sterility, after a limited period of prosperity.

The development of mathematics is realized not only through mathematical proofs, but also through the axioms and definitions, which, one might say, are arbitrary to a certain extent. It is here that we clearly see the synthetic character (axioms) and the formalism of the mathematics method. One of the obvious examples are the non-Euclidean geometries, having risen from considerations within geometry related to formally proving Euclidean fifth postulate. Very abstract discussions for that time (and for many today too) found unanticipated applications in mathematics, physics and other human activities, thus changing fundamentally the world picture known thus far. Will geometry of Lobachevsky become commercialized one day and as such the source of adoration of the laity and the infinite satisfaction of the merchants? With the affirmative answer, we add that the one hundred years distance has been stretched and that we will not wait too long for the mass consumption.

To be abstract means to overcome special, concrete cases of actual happening. But to overcome an actual case, does not mean to completely sever ties with it. It is important to restrict oneself to a certain set of abstractions at a given moment. The advantage of restricting one's attention to such a set of abstractions is in restricting our thinking to completely and clearly specific things, with completely and clearly determined relationships. Accordingly, if one is logically oriented, one may deduce a sequence of various conclusions concerning the relationships among these abstract entities. Moreover, if those abstractions are well founded, that is, if they do not abstract everything important in our experience, the scientific thought restricted to them will arrive to a sequence of various important truths in connection to our perception of nature.

There are many misunderstandings related to the abstractness of mathematics, such as that it is "useless," etc., often because of the inability to understand and adopt those abstractions and then because a long road needs to be passed for the abstractions to be perceived by the "senses." Much can be done for resolving this problem: It is known how Steinhaus formulated the famous and extremely abstract theorem of Brouwer through the statement that a sphere cannot be ideally combed. Or another such a statement, known as the statement that every sandwich may be cut with one straight cut, in a way that bread, cheese and salami are cut into two equal pieces. Although these statements are popularized in this way, they are at the same time too oversimplified, and it is a fact that one should not exaggerate in approaches like this to mathematics. Yet another example is representation of a five-dimensional space: if every spatial coordinate is accompanied by time, we already get the four-dimensional space; pressure in every such a point in space and time may represent the fifth dimension. And so we can go on, using additional physical and sensual characteristics as additional dimensions.

Mathematics, however, could not survive if it lived only through the sensory observations. Otherwise, what could be said about infinity, about non-Euclidean geometries, existence theorems, etc?

Some consider important to be only those problems and mathematics results that can somehow be useful to physics, thus indirectly to the sensory existence. However, in principle any mathematics theorem may find applications in physics, for every form of logical reasoning may prove to be useful in investigations of the real world. For instance, those who consider the theory of differential equations (or today more modern results of differential geometry) to cater more for the needs of mechanics than some other theories, judge this from the point of view of the modern state of mathematics and mechanics. There is no serious argument to support the claim that some areas of mathematics are better suited than the others in pursuit of studying nature. Although the theories of ordinal numbers or the theory of projective dimension of modules do not find applications in physics, it just means that the physics theories are now formulated in the form that does not allow application of ordinal numbers to them. On the other hand, Statistical Mechanics (or Financial Engineering for that matter) uses at least one of the fresh achievements of set theory – Lebesgue measure, and crystallography and relativity theory use very extensively – group theory. The more ready-made mathematical forms and methods we have at our disposal, the easier it will be to a physicist to choose what is needed for his work, at the given moment. In this sense, the needs of science are similar to ordinary human needs: The very fact that there are many different means for their satisfaction, serves as a stimulus of the further growth of the needs.

The development of mathematics goes most frequently along the road of abstraction open only to great minds who are able to peek in a perfectly different way in an approach to problems posed before mathematics - resulting in a creation of the whole mathematical theories. We state several examples of such great discoveries: Euclid's *Elements*, Archimedes' foundations of integral calculus, Descartes' analytic geometry, Leibnitz-Newton differential calculus, non-Euclidean geometry of Lobachevsky, Bolyai and Gauss, Cauchy's theory of complex variables, theory of Riemann surfaces and his paper "On hypotheses that lie in the foundations of geometry", group theory of Abel and Galois, Cantor's set theory, modern topology, category theory of Eilenberg and Mac Lane. The history of mathematics does not contain many great jumps like these. Their significance is in that after each of them, there comes a century of concentrated work of other mathematicians who systematize and derive consequences from the theories created by the genii.

Apart from a very characteristic and developed mathematical language by which one can express, in a short form, the same or more than other languages express through hundreds of filled pages (and which must be studied for a long time too, in order to master it completely), the difficulty in understanding mathematics literature is also in reasoning, which the authors of the papers never present completely, but only in a (sometimes drastically) shortened form, for the full proofs would enlarge the volumes of published works several tens of times. The volume of mathematics literature and multitude of the references in every paper only enlarge the problem of unhampered understanding of mathematical results. There were only a couple of periodic mathematics journal in the XVIIth century (*Journal des sçavans* and *Philosophical Transactions of the Royal Society*), there were around 600 of them around 1920, and there are several thousands of them today and the new arrive daily to mathematics libraries. This is at the same time an answer to layman's question whether mathematics is developing (for they say: "Everything in mathematics has already been known since the times of Euclid!"). Mathematics must also face the fact that among our contemporaries there are people whose mathematical knowledge belongs to the epoch older than the Egyptian pyramids and they are greatest in number. Mathematical knowledge of a small number of people reaches to the epoch of the middle ages, and the level of mathematics of the XVIIIth century is not attained even by one in thousands. We have already mentioned that if we want to transform a primitive man into a mathematician, we cannot rely on evolution, because the period that separates primitive from modern man is too short. The human brain cannot be transformed by a move of a magic baton. In order to make Poincaré from a cave man, we must educate generation after generation, along the thorny road, that it seems cannot be shortened. However the distance between those who belong to the avant-garde and the countless mass of the passengers is growing more and more, the chasm is bigger and bigger and those who are leading are further and further away. They are disappearing from sight, they are not known by many, unusual stories are told about them. There are even those who do not believe in their existence.

At the end of this section we again quote Hardy [5] in relation to the judgment of mathematics (he pursued): "The 'real' mathematics of the 'real' mathematicians, the mathematics of Fermat and Euler and Gauss and Abel and Riemann, is almost wholly 'useless' (and this is as true of 'applied' as of 'pure' mathematics). It is not possible to justify the life of any genuine professional mathematician on the ground of the 'utility' of his work."

Before we go into discussions on mathematics as a reality, we note that Hardy too was sometimes very practical. For instance his scale of mathematicians is well known – it is formed according to the difficulty of a statement a mathematician in question proved or discovered. His weakest student was given 1 and Einstein was given 100.

## Mathematics: Reality and the sensory world

We will say here a few words about the relationship of mathematics and the physical world in order not to fall into the situation like that when Honoré de Balzac who, interrupting his conversation, said: "It is fine my dear friend, but let's come back to reality, let's talk about Eugénie Grandet..."

It can happen that we understand all the abstract facts about the Sun, and everything about the atmosphere and everything about the rotation of the Earth, and still miss the splendor of a sunset. A direct observation of a concrete phenomenon in its actuality cannot be replaced by anything. Just as instinct without feeling constitutes crudeness, feelings without instinct constitute decadence.

It would be a great mistake to think that any given invention consists of a bare scientific idea, that just has to be picked up and used. From an idea to its application there is a period of vivid imaginative projecting. One of the elements in a new method is exactly in finding the ways to start bridging the split between a scientific idea and a final realization. It is a process of disciplined overcoming of difficulties, one after another.

A great part of the ancient mathematics (still taught today in schools and universities) has come about through a practical, everyday world. Geometry (nowadays called Euclidean geometry) has come about as is well known, and as its own name says, from the needs of everyday land measurements (in a difference from non-Euclidean geometries, that did not have a previous anticipation in human experience). Another example is Archimedes' anticipation of integral calculus while solving the problem of volumes of barrels.

Thus it is possible to model physical phenomena by mathematical models (which is not the least easy) and arrive to new mathematical theories and new problems. When we model and describe, in mathematical language, movement of a pendulum of length $l$, under the Earth's gravity $g$ (clearly, under ideal conditions, without resistance of air for instance) and the angle $x(t)$ of the deviation (at time $t$) from the vertical position, then we arrive at the following differential equation:

$$l\, x'' + g \sin x = 0$$

which cannot be solved in an elementary way (with elementary functions); the attempts to solve this equation lead to a new theory of elliptic integrals. Or consider the concrete example in cartography: To map some area *A* onto a given region *B* in a way that small parts of *A* are mapped into small parts of *B* and that the boundary of one region goes into the boundary of another. Mathematical modeling of this request leads to an abstract theorem of Riemann on conformal mappings, which continues to exist further, independently, as a mathematics theory, and moreover it can be used again in treatment of other practical problems.

If we consider simple problems known to everybody (with easily understood formulations) as simple reality, then it is certain that in generalizing the Pythagoras' theorem to Fermat's last problem, we acted very fruitfully for the attempts to solve the latter created important mathematical theories.

Mathematical theories that arose in describing some natural phenomena could be in return used later for description of other phenomena. Thus, for instance, the differential equation used in modeling water waives was used for a description of acoustic and radio waives. The stimulus for Cantor's set theory was in seeking answers to some problems of Fourier series, which in turn arose from practical considerations. Probability theory arose at the beginning as a purely practical discipline connected to the possibilities in gambling games, but today, highly abstract, it has considerable successes in explaining concrete physical and societal phenomena.

Mathematics is a miraculous creation, and the miracle consists in the fact that the human brain is able to construct models of complex natural phenomena, which, it may at first seem do not yield to a description in mathematical terms. Natural sciences begun with mathematics and would not continue to exist long after throwing mathematics out of their content (if such an exception would be possible at all). Our century has seen a proliferation of laboratories for a mass production of data. Whether the obtained data will remain only data or they would transform into science, will depend on the degree of interactions with the spirit of mathematics they will achieve. We may and should use specific mathematical methods (not seldom incorrectly called the methods of applied mathematics) in places where there is talk about applications of practical results, demonstration of calculations or the use of formulas derived with the aid of purely logical reasoning.

It is quite probable that if mathematics feeds on its own juices, then the stimuli that drive it will considerably diminish, and possibly completely exhaust in time.

We arrive inevitably to the question of "applied" and "pure" mathematics and "applied" and "pure" mathematicians. It is clear that this question is academic in its essence (but not only academic), for, instead of distinguishing these terms, we should make a difference between mathematics without content and mathematics with content, mathematics with the goal and aimless mathematics, and then as Goethe said, "a men has freedom to pursue that which most attracts him, which gives him satisfaction, which appears to him to be most useful."

Often, "applied" mathematics is identified with the application of some mathematical construction to physical world. If we accept this, we must also accept that a waiter, when he brings us the bill (with his mathematics, more or less well applied) is an applied mathematician. Likewise, am I an applied chemist, if I am writing this by a ball-point pen? We understand the term "applied" mathematics here in a creative sense, which results in its non-existence as such, because it falls under the general terms "mathematics" and "mathematician." If a natural phenomenon is described by a mathematical apparatus, it is most often an already known apparatus, while we can talk about creative mathematics (thus also about "applied") only if such a modeling breeds independent and new abstract mathematical model that exists independently, regardless of the physical world from which it arose and from which new mathematical use (arrived at inside the mathematical apparatus) may be drawn. The price of such a mathematical theory is greater if, when returning to "earth" it can describe other, up to that point non-unified phenomena, apart from the phenomenon it arose from. Thus we can talk about the unique mathematics, but also that mathematics realizes a synthesis with other scientific disciplines.

We give an example of Milutin Milankovitch, a Serbian scientist who gave extraordinary and thus far unsurpassed results in climatology and theory of ice ages [10, 11]; those results are very appreciated and much used in today's development of these sciences. In studying climate and the ice ages Milankovitch considers specially the influence of astronomical and geophysical elements on climate. The apparatus he uses in his work shows that we have a scientist who is well informed in mathematics, physics, astronomy and mechanics. The mathematics he uses encompasses the area nowadays studied in university courses of Analysis I and II and it was enough for Milankovitch to determine by mathematical tools the quantity of heat that arrives to some celestial body, depending on various variable quantities connected with that body. Milankovitch did not make an abstract mathematical theory out of problems he treated, which would live in mathematics independently of those problems. Also, the theory he developed has not been used for modeling of other physical and societal phenomena, thus we cannot say that he was a mathematician (neither an "applied" mathematician). In the synthesis he realized, Milankovitch can be most suitably defined as a mathematical climatologist (whatever this term means), who used mathematical methods in his investigations, applied in a new ambient of climatology. We also think that a "mathematical biologist," or, for instance, "mathematical physicist" are more suitable terms than an "applied mathematician." It is our belief that Milankovitch' theory will in the future be an anticipation of a mathematical theory which will exist abstractly and whose results would be applied in mathematics and the sensory world. The contemporary works in catastrophe theory (Gilmore, [4]) undoubtedly go in the direction of creation of such a theory.

It is a dilemma today how in ever increasing specialization a mathematician would get introduced into problems of physical world and its modeling. It would be necessary in the future to have more and more theoretical physicists, biologists and others who would acquire knowledge about deep cognition of mathematical principles, as well as mathematicians who will not restrict themselves only to the aesthetic development of mathematics abstractions. Collaboration in examining problems between mathematicians and experts in other disciplines does not start with asking the problem, but much earlier, at the initial stage of the idea and the joint formulation of the problem. An approximate quantity needed for that is given by the Steinhaus' number of 50 hours of common work on establishing the problem and its mathematical formulation [14].

We emphasize that systematic creativity can occur in every field, under the conditions where talent and high professionalism exist. And these are certainly improved through love towards work which makes us see and highly esteem the beauty i. e. harmony and functionality, which are fundamental categories of all creative activity and especially mathematics, which stands firmly on earth with its head high in the clouds.

# References

The references are not systematic and serve to direct the reader to information about the themes we have touched upon here to some extent.